# Elaboration on Two Points Raised in "Classifier Technology and the Illusion of Progress"

**Robert C. Holte**

## 1. INTRODUCTION

This short note elaborates two points raised in David Hand's target article. First, I provide additional evidence that simple classification rules should be given serious consideration in any application and that there are often diminishing returns in considering increasingly complex classifiers. Second, I refine Hand's basic argument that small improvements in performance are irrelevant because of the uncertainty about many aspects of the situation in which the classifier will be deployed. In particular, I briefly describe a recently developed method for analyzing and comparing classifier performance when the class ratios and misclassification costs are unknown. This does not refute his general argument, but it does provide an important exception to it.

## 2. SIMPLICITY-FIRST METHODOLOGY AND DIMINISHING RETURNS

Hand (Section 2.3) cites my 1993 study [4] in which the accuracy of one-level decision trees, which classify examples based on the value of a single feature, was compared to the accuracy of the decision trees learned by C4.5 [8], a state-of-the-art decision tree learning algorithm. The article caused quite a stir, because nobody at the time suspected that most of C4.5's classification accuracy could be achieved, on many of the standard test data sets, by building just the first level of the decision tree. The overall conclusion of my 1993 article is the same as Hand's—not that the more complex decision rules should be cast aside, but that the simple decision rules should not be dismissed out of hand. One can never tell, a priori, how much of the structure in a domain can be captured by a very simple decision rule, and since simplicity is advantageous for both theoretical and practical reasons, it is incumbent on a responsible experimentalist or practitioner to begin with the simplest decision rules. Only if they prove unacceptable should more complex decision rules be considered. I coined the term "simplicity-first methodology" to describe this systematic approach of proceeding from simple to more complex decision rules.

In a follow-up paper [1], Maass and Auer developed an efficient algorithm for constructing a decision tree of fixed depth $d$, with the minimal error rate on the training data, and we proved theoretical bounds on the generalization error rate of this decision tree. This empirical study showed that the performance advantage of C4.5 over one-level trees in my original study [4] greatly diminishes when depth is increased to two, with the two-level trees actually being superior to C4.5's trees on 4 of the 15 data sets in the study.

Table 1 herein compares the accuracies achieved when $d = 0$, $d = 1$ and $d = 2$. These accuracies are averages of nine repetitions of 25-fold cross-validation on each data set. The $\Delta(1\text{–}0)$ column gives the accuracy improvement achieved by moving from a zero-level tree, which classifies all examples according to the majority class, to a one-level tree, and the $\Delta(2\text{–}1)$ column gives the accuracy improvement achieved by moving from a one-level tree to a two-level tree. Comparing these two columns, we see clear confirmation of Hand's observation that increasing complexity produces diminishing returns on accuracy improvement in many domains.

There have been other studies that showed that simple classifiers perform well on standard test data sets. Domingos and Pazzani [2] showed that a naive Bayesian classification algorithm significantly outperformed state-of-the-art systems for decision tree learning, decision rule learning and instance-based learning in a substantial number of the 28 data sets in their study. Kohavi [5] showed that wrapper-based feature selection, combined with a majority classifier, can produce simple classifiers that are as accurate as C4.5's trees in many cases. Linear discriminants (perceptrons) have also been seen to perform surprisingly well [6, 9].

*Robert C. Holte is Professor, Department of Computing Science, University of Alberta, Edmonton, Alberta, Canada, T6G 2E8 (e-mail: holte@cs.ualberta.ca).*

TABLE 1
*Diminishing returns with additional complexity*[*]

| Data set | Zero-level | One-level | Two-level | $\Delta$(1–0) | $\Delta$(2–1) |
|---|---|---|---|---|---|
| BC | 70.3 | 67.2 | 66.3 | −3.1 | −0.9 |
| HE | 79.4 | 79.2 | 78.6 | −0.2 | −0.6 |
| AP | 80.2 | 80.0 | 88.6 | −0.2 | 8.6 |
| SE | 90.7 | 95.0 | 97.3 | 4.3 | 2.3 |
| LA | 64.9 | 71.6 | 86.6 | 6.7 | 15.0 |
| PI | 65.1 | 73.6 | 74.8 | 8.5 | 1.2 |
| SP (3) | 51.9 | 63.2 | 79.4 | 11.3 | 16.2 |
| CH | 52.2 | 66.1 | 86.9 | 13.9 | 20.8 |
| IO | 64.1 | 78.3 | 86.1 | 14.2 | 7.8 |
| PR | 50.0 | 66.3 | 69.3 | 16.3 | 3.0 |
| HD | 54.5 | 70.9 | 67.1 | 16.4 | −3.8 |
| G2 | 53.4 | 76.2 | 79.7 | 22.8 | 3.5 |
| CR | 55.5 | 85.5 | 84.2 | 30.0 | −1.3 |
| SO (4) | 36.2 | 85.3 | 91.1 | 49.1 | 5.8 |
| IR (3) | 33.3 | 91.9 | 95.7 | 58.6 | 3.8 |

[*]The first column gives the acronym for the data set as in [1], with the number of classes shown in parentheses if it is different from two. The next three columns give the accuracy of the majority classifier (zero-level decision tree), one-level decision tree and two-level decision tree, respectively. The $\Delta$(1–0) column gives the difference in accuracy between the one-level and zero-level trees, and the final column gives the difference in accuracy between the two-level and one-level trees. The rows are sorted according to $\Delta$(1–0).

## 3. EMPIRICAL COMPARISONS OF CLASSIFIERS IN UNKNOWN CIRCUMSTANCES

The fundamental argument put forward by David Hand has two parts: (1) that often only small performance gains arise from using complex classifiers and (2) that the small gains seen in the idealized laboratory setting will be swamped, in practical applications, by unpredictable and changing conditions that have a substantial effect on performance. I agree with both of these statements, in general, but I would like to point out, with regard to the latter, that we do possess methods for coping perfectly well with certain important kinds of unpredictable and changing circumstances.

Among the most important examples Hand gives of unpredictable and changing factors that affect a classifier's usefulness in practice are the costs of the different types of misclassification and the distribution of data to which the classifier will be applied. I agree entirely that in many practical settings these factors cannot be determined at the time classifiers are being evaluated and compared, and that these factors often change with time.

Drummond and I have developed a method, called cost curves, for analyzing and comparing two-class classifier performance when the misclassification costs and the relative frequency of the two classes are unknown [3]. The key idea is to plot performance (expected cost, normalized to be between 0 and 1) as a function of these unknowns. It turns out that, for the case of expected cost, these unknowns can be combined into a single aggregate unknown that also varies between 0 and 1. Cost curves therefore are a two-dimensional plot, with performance (normalized expected cost) as the $y$-axis and the aggregate unknown, which we call $PC(+)$, as the $x$-axis.

The cost curve for a given classifier is a straight line that depicts its performance across all possible combinations of misclassification costs and class ratios. Empirical confidence intervals can be computed for cost curves and for differences between cost curves, allowing one to answer the all-important question, "Under what circumstances does classifier $A$ significantly outperform classifier $B$?" A software tool that fully supports cost curve analysis is available upon request.

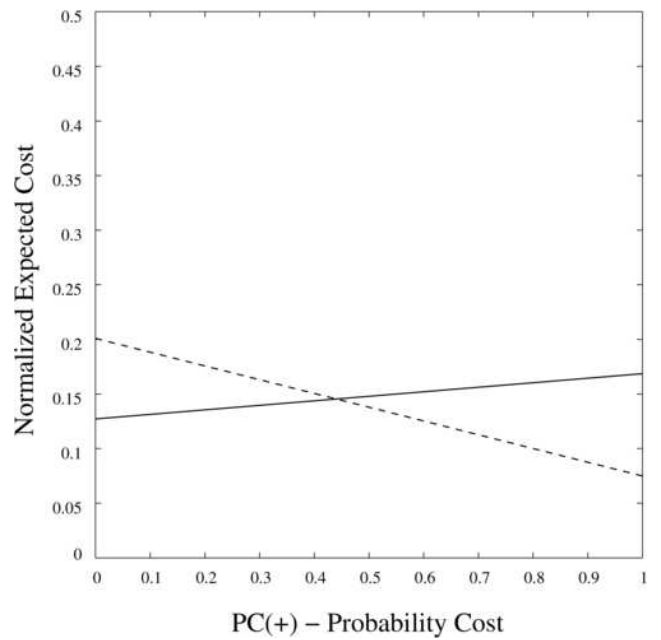

FIG. 1. *Cost curves for* C4.5 *(solid line) and* 1R *(dashed line) on the Japanese credit screening data set.*

Figure 1 herein shows the cost curves for two classifiers on the Japanese credit screening data from the UCI repository [7]. The solid line is the cost curve for C4.5's decision tree on this data set and the dashed line is the cost curve for the one-level decision tree produced by my $1R$ system [4]. We can see that these two classifiers have identical performance when $PC(+)$ has a value of roughly 0.45, that the one-level tree has a lower expected cost than C4.5's decision tree for larger values of $PC(+)$ and that C4.5's tree outperforms the one-level tree for smaller values.

My aim here is not to give a tutorial on cost curves, but to point out that there are sound, practical ways to cope with some of the factors that Hand correctly identifies as often being unknown, or subject to change, at the time of classifier evaluation. Cost curves provide a concrete example of how we can do classifier evaluation and comparison perfectly well without any knowledge about misclassification costs or the class ratios. By considering all possible combinations of the unknown factors, exact analysis and comparison is possible, and small performance differences can be significant. However, this does not refute Hand's general point. There are other factors and kinds of changes, such as shifting distributions within a class [10], that we do not yet know how to cope with—a challenge for future research.

## ACKNOWLEDGMENTS

I would like to thank the Natural Sciences and Engineering Research Council of Canada for financial support, and Alberta Ingenuity for funding the Alberta Ingenuity Centre for Machine Learning. All my work on cost curves is joint with Chris Drummond of the Institute for Information Technology (Ottawa) of the Canadian National Research Council.